\documentclass[12pt]{amsart}
\usepackage[usenames,dvipsnames]{xcolor}
\usepackage{tikz}
\usepackage{fancyhdr}
\usepackage{txfonts}
\usepackage{graphicx}
\usepackage{epsfig}
\usepackage{mathrsfs}
\usepackage{amssymb}
\usepackage{latexsym}
\usepackage{amsmath}
\usepackage{amsfonts}
\usepackage{amsbsy}
\usepackage{amscd}
\usepackage{subfigure}
\usepackage{verbatim}
\usepackage{color,tikz}
\usetikzlibrary{calc}

\usepackage{amsfonts}
\usepackage{amsbsy,calc}
 \usepackage{graphicx,ifthen,cite}
\usepackage{dsfont}

\newcommand{\A}{{\mathcal A}}
\newcommand{\B}{{\mathcal B}}
\newcommand{\D}{{\mathcal D}}

\newcommand{\ZZ}{{\mathcal Z}}
\newcommand{\R}{{\mathbb R}}
\newcommand{\Z}{{\mathbb Z}}
\newcommand{\C}{{\mathcal C}}

\newcommand{\de}{\delta}

\newcommand{\la}{\lambda}
\newcommand{\La}{\Lambda}

\newtheorem{prop}{Proposition}[section]
\newtheorem{lem}[prop]{Lemma}

\newtheorem{defi}[prop]{Definition}
\newtheorem{coro}[prop]{Corollary}
\newtheorem{theo}[prop]{Theorem}

\newcommand{\dev}[2][]{{\,\textup{d}}^{#1}#2}

\def\dmu{\dev\mu}

\newif\ifdraft
\drafttrue
\ifdraft
\else
\fi
\numberwithin{equation}{section}
\begin{document}
\title[The generalized Fuglede's conjecture holds for a class of Cantor-Moran measures]{The generalized Fuglede's conjecture holds for a class of Cantor-Moran measures}
\author[L-X An, Q Li, M-M Zhang]{Li-Xiang An, Qian Li$^{\ast}$, Min-Min Zhang}
\address{School of Mathematics and Statistics, $\&$ Hubei Key Laboratory of Mathematical Sciences, Central China Normal University, Wuhan 430079, P.R. China.}
\email{anlixianghai@163.com}
\address{School of Mathematics and Statistics, Huazhong University of Science and Technology, Wuhan 430074, P. R. China.}
\email{liqian303606@163.com}
\address{School of Mathematics and Statistics, $\&$ Hubei Key Laboratory of Mathematical Sciences, Central China Normal University, Wuhan 430079, P.R. China.}
\email{zhangminmin0907@163.com}

\thanks{*Corresponding author.}
\thanks{ This work was supported by the National Natural Science Foundation of China (No. 12171181 and 12371087).}

\maketitle
\begin{abstract}
Suppose ${\bf b}=\{b_n\}_{n=1}^{\infty}$ is a sequence of integers bigger than 1 and ${\bf D}=\{{\mathcal D}_{n}\}_{n=1}^{\infty}$ is a sequence of consecutive digit sets. Let $\mu_{{\bf b},{\bf D}}$ be the Cantor-Moran measure defined by
\begin{eqnarray*}
  \mu_{{\bf b},{\bf D}}&=& \delta_{\frac{1}{b_1}\D_{1}}\ast\delta_{\frac{1}{b_1b_2}\D_{2}}\ast
\delta_{\frac{1}{b_1b_2b_3}\D_{3}}\ast\cdots.
  \end{eqnarray*}
   We prove that $L^2(\mu_{{\bf b},{\bf D}})$ possesses an exponential orthonormal basis if and only if $\mu_{{\bf b},{\bf D}}\ast\nu={\mathcal L}_{[0,N_1/b_1]}$ for some Borel probability measure $\nu$.
  This theorem shows that  the generalized Fuglede's conjecture is true for such Cantor-Moran measure. An immediate consequence of this result is the equivalence between the existence of an exponential orthonormal basis and the integral tiling of ${\bf D}_n=\D_{n}+b_n\D_{n-1}+b_2\cdots b_n\D_{1}$ for $n\geq1$.
\end{abstract}

\section{Introduction}

A  Borel probability measure $\mu$ on ${\mathbb R}^d$ is called a {\it spectral measure} if we can find a countable set $\Lambda\subset{\mathbb R}^d$ (called {\it spectrum}) such that the set of exponential functions $E(\Lambda): = \{e^{2\pi i \lambda \cdot x}:\lambda\in\Lambda\}$ forms an orthonormal basis for $L^2(\mu)$.  If $\Omega\subset\R^d$ is a measurable set with finite positive Lebesgue measure and ${\mathcal L}_{\Omega}$ is a spectral measure, then we say that $\Omega$ is a {\it spectral set}. Here ${\mathcal L}_{K}$ denotes the normalized Lebesgue measure restricted to the measurable set $K$ of finite positive Lebesgue measure.

It is well-known from classical Fourier analysis that the unit cube $[0,1]^d$ is a spectral set with a spectrum ${\mathbb Z}^d$. Which other sets $\Omega$ can be spectral? The research on this problem has been influenced for many years by a famous paper \cite{Fu1974} due to Fuglede, who suggested that there should be a concrete, geometric way to characterize the spectral sets.

{\bf Conjecture (Fuglede's Conjecture):}
A set $\Omega\subset\R^d$ is spectral if and only if it can tile the space by translations.

We say $\Omega$ is a translational tile if there exists a discrete set ${\mathcal J}$ such that the translated copies $\left\{\Omega+t: t\in{\mathcal J}\right\}$ constitute a partition of $\R^d$ up to measure zero. Although the conjecture was disproved eventually in dimension 3 or higher in
its full generality\cite{KM20061, KM20062, T2004}, most of the known examples of spectral sets are constructed from translational
tiles. An important result is proved by Lev and  Matolcsi \cite{LM2022} that all
spectral sets must admit a ``weak tiling" which is a generalization of translational tiling in its measure theoretic form. Let $K\subset\R^d$ be a  bounded, measurable set. We say that another measurable, possibly unbounded, set $\Sigma\subset\R^d$ admits a weak tiling by translates of $K$, if there exists a positive, locally finite (Borel) measure $\nu$ on $\R^d$ such that ${\bf 1}_K\ast\nu={\bf 1}_{\Sigma}$, where ${\bf 1}_{A}$ denotes the indicator function of a set $A$.

Jorgensen and Pedersen \cite{JP1998} widened the scope of Fuglede's conjecture and  discovered that the standard middle-fourth Cantor measure $\mu_{1/2k,\{0,2\}}$ is a spectral measure. It is the first spectral measure that is non-atomic and singular measure.
Strichartz \cite{Str2000, Str2006} discovered a surprising and interesting phenomenon: the Fourier series corresponding to certain spectra of $\mu_{4^{-1},\{0,2\}}$ can exhibit significantly better convergence properties than their classical counterparts on the unit interval. Specifically, the Fourier series of continuous functions converge uniformly and Fourier series of $L^p$-functions converge in the $L^p$-norm for $1\le p<\infty$.
Following these discoveries, there has been considerable research on such measures \cite{AW2021,DH2016,LW2002,FHL2015,DHS2009,D2012,DHL2013,DHL2014,Dai2016,HL2008,DHL2019}, and a celebrated open problem was to characterize the spectral property of the Cantor measures $\mu_{\rho,\D}, 0<\rho<1$  among the $N-$Bernoulli convolutions
$$\mu_{\rho, \D}(\cdot)=\frac{1}{N}\sum_{d\in\D}\mu_{\rho,\D}(\rho^{-1}(\cdot)\cdot-d),$$
where $\D=\{0,1,\ldots,N-1\}$. Hu and Lau \cite{HL2008}, as well as  Dai \cite{D2012} showed that the above Cantor measures $\mu_{1/(2k),\{0, 2\}}$ are the only class of spectral measures among the $\mu_{\rho,\{0,2\}}$. Dai, He and Lau \cite{DHL2014} demonstrated that a similar result holds for $N-$Bernoulli measures $\mu_{\rho, \D}$.



Later on, Strichartz \cite{Str2000} has already formulated the most general  fractal spectral measures one can possibly generate. Let ${\bf b}=\{b_n\}_{n=1}^{\infty}$ be a sequence of integers bigger than 1 and let $\textbf{D}=\{\D_{n}\}_{n=1}^{\infty}$ be a sequence of integer digit sets.
Let $\delta_a$ be the Dirac measure and
denote
$$\delta_E=\frac 1{\#E}\sum_{e\in E}\delta_e$$
for a finite set $E$. Write
 $$
\mu_{n} = \delta_{\frac{1}{b_1}\D_1}\ast\delta_{\frac{1}{b_1b_2}\D_2}\ast
\delta_{\frac{1}{b_1b_2b_3}\D_3}\ast\cdots\ast\delta_{\frac{1}{b_1b_2\cdots b_n}\D_n},
$$
where $\ast$ is the convolution sign. If the sequence of
convolutions $\{\mu_n\}_{n=1}^{\infty}$ converges weakly to a Borel probability measure $\mu_{{\bf b},\textbf{D}}$ with compact support, then we call $\mu_{{\bf b},\textbf{D}}$ a {\it Cantor-Moran measure} as a generalization of the standard Cantor measure studied first by Moran \cite{M1946}. {This opens a research field for orthogonal harmonic analysis of Cantor-Moran measures (e.g. \cite{AFL2019,AH2014,DL2022,LMW2022}).} In this paper we study the spectrality of the Cantor-Moran measures generated by
an integer sequence ${\bf b}=\{b_n\}_{n=1}^{\infty}$ with $b_n\geq2$ and a sequence of consecutive digit sets $\textbf{D}=\{\D_{n}\}_{n=1}^{\infty}$, i.e., $\D_{n}=\{0,1,\ldots, N_n-1\}$. We first provide a sufficient and necessary condition for the existence of such Cantor-Moran measure.

\begin{theo}\label{th1.1}
Let ${\bf b}=\{b_n\}_{n=1}^{\infty}$ be a sequence of integers bigger than 1 and  ${\bf D}=\{\D_{n}\}_{n=1}^{\infty}$ be a sequence of consecutive digit sets with $\D_{n}=\{0,1,\ldots, N_n-1\}$ where $N_n\geq2$. Then  the sequence of discrete measures $$\mu_n=\delta_{\frac{1}{b_1}\D_{1}}\ast\delta_{\frac{1}{b_1b_2}\D_{2}}\ast\cdots\ast
\delta_{\frac{1}{b_1b_2\cdots b_n}\D_{n}}$$ converges weakly to a Borel probability measure $\mu_{{\bf b}, {\bf D}}$ if and only if
\begin{equation}\label{eq1.1}
  \sum_{n=1}^{\infty}\frac{N_n}{b_1b_2\cdots b_n}<\infty.
\end{equation}
In this case,  $\mu_{{\bf b}, {\bf D}}$  supports on a compact set
$$T({\bf b}, {\bf D})=\left\{\sum_{n=1}^{\infty}\frac{d_n}{b_1b_2\cdots b_n}:\,d_n\in\D_{n}\right\}:=\sum_{n=1}^{\infty}\frac{\D_{n}}{b_1b_2\cdots b_n}.$$
\end{theo}

The spectral property of such measure was first studied by the first-named author and He  \cite{AH2014} as a generalization of the $N-$Bernoulli convolution (all $b_n=b, D_n=\{0,1,\ldots,N-1\}$). The first-named author and He showed that  $\mu_{{\bf b}, {\bf D}}$ is spectral when $N_n$ divides $b_n$ for each $n\geq1$. Under the condition that $\{N_n\}_{n=1}^{\infty}$ is bounded, it has been proved in \cite{DL2022} that $N_n$ divides $b_n$ for each $n\geq1$ is also a necessary conditions for $\mu_{{\bf b}, {\bf D}}$ to be spectral.
These Cantor-Moran measures also show that spectral measures can have support of any Hausdorff dimensions \cite{DS2015}. Last but not least, Cantor-Moran measures offer new examples of fractal measures that admit Fourier frame but not a Fourier orthonormal basis \cite{GL2014}, which lead to a new avenue to study a long-standing problem whether a middle-third Cantor measure has a Fourier frame.

Our motivation to extend the $N-$Bernoulli convolutions to this class of measures due to the conjecture by Gabardo and Lai \cite{GL2014} and aims to answer a question about the Cantor-Moran spectral measure and the integer tiles.  To describe a unifying framework bridging the gap between singular spectral measures and spectral sets, Gabardo and Lai \cite{GL2014} extended the classical Fuglede's Conjecture to a more generalized form.


{\bf Conjecture (Generalized Fuglede's Conjecture):} A compactly supported Borel
probability measure $\mu$ on $\R$ is spectral if and only if there exists a Borel probability $\nu$ on $\R$ and a fundamental domain $Q$ of some lattice on $\R$ such that $\mu\ast\nu={\mathcal L}_{Q}$.


Deterministic positive results about Cantor measure have been appeared in many papers (e.g. \cite{AW2021,DHS2009,D2012,DHL2013,DHL2014,Dai2016,HL2008,DHL2019,LW2002}). However, there are relatively few results regarding Cantor-Moran measures. Gabardo and Lai \cite{GL2014} showed that if both $\mu$ and $\nu$ are two singular probability measures with $\mu\ast\nu={\mathcal L}_{[0, 1]}$, then they all are Cantor-Moran measures studied by us.
In this paper, we will demonstrate that the generalized Fuglede's Conjecture holds for our target measure.

\begin{theo}\label{th1.2}
Suppose ${\bf b}=\{b_n\}_{n=1}^{\infty}$ is a sequence of integers bigger than 1 and  ${\bf D}=\{\D_{n}\}_{n=1}^{\infty}$ is a sequence of consecutive digit sets with $\D_{n}=\{0,1,\ldots, N_n-1\}$ where $N_n\geq2$. Then the following are equivalent.

{\rm(i)} The Cantor-Moran measure $\mu_{{\bf b}, {\bf D}}$ is spectral.

{\rm(ii)} There exists a Borel probability $\nu$ such that $\mu_{{\bf b}, {\bf D}}\ast\nu={\mathcal L}_{[0,N_1/b_1]}$.

{\rm(iii)} $N_n$ divides $b_n$ for each $n\ge2$.
\end{theo}

In fact, the implication $``{\rm(ii)}\Rightarrow {\rm(i)}"$ stems from \cite[Theorem 1.1]{GL2014 }, and the proof of $``{\rm(iii)}\Rightarrow {\rm(i)}"$ is provided in\cite[Theorem 1.4]{AH2014}, while $``{\rm(i)}\Rightarrow {\rm(iii)}"$ and $``{\rm(iii)}\Rightarrow {\rm(ii)}"$ are apparently new in this generality. We now outline the strategy of the proof. Let us set up the notations.
\begin{eqnarray*}
\mu_{{\bf b}, {\bf D}}&=&\delta_{\frac{1}{b_1}\D_{1}}\ast\delta_{\frac{1}{b_1b_2}\D_{2}}\ast
\delta_{\frac{1}{b_1b_2b_3}\D_{3}}\ast\cdots\\
&=&\mu_{n}\ast\mu_{>n},
\end{eqnarray*}
where $\mu_n$ is the convolutional product of the first $n$ discrete measures and $\mu_{>n}$ is the
remaining part. Based on the aforementioned decomposition of $\mu_{{\bf b}, {\bf D}}$, one seeks to observe the behavior of $\mu_n$ and $\mu_{>n}$ under the condition that $\mu_{{\bf b}, {\bf D}}$ is a spectral measure. Actually, if $0\in\Lambda$ is a spectrum of ${\mu}_{{\bf b}, {\bf D}}$ and for any $n\ge1$, then we can construct spectra of $\mu_n$ and $\mu_{>n}$ relying on the suitable decomposition of $\Lambda$(see Definition \ref{defi1} for detail).

\begin{theo}\label{th1.3}
Suppose ${\bf b}=\{b_n\}_{n=1}^{\infty}$ is a sequence of integers bigger than 1 and  ${\bf D}=\{\D_{n}\}_{n=1}^{\infty}$ is a sequence of consecutive digit sets with $\D_{n}=\{0,1,\ldots, N_n-1\}$ where $N_n\geq2$. If $0\in\Lambda$ is a spectrum of ${\mu}_{{\bf b}, {\bf D}}$, then for each $n\geq1$, we can decompose it suitable, denoted by $\Lambda=\bigcup_{\alpha\in \A}\Lambda_\alpha$, such that $\A$ is a spectrum of $\mu_n$ and each $\Lambda_{\alpha}$ is a spectrum of $\mu_{>n}$.
\end{theo}

Applying Theorem \ref{th1.3}, we can reduce $``{\rm(i)}\Rightarrow {\rm(iii)}"$ to the following theorem, considering $\mu_n$ as a spectral measure.


\begin{theo}\label{th1.4}
The discrete measure $\mu_n=\ast_{j=1}^n\delta_{\frac1{b_1\cdots b_j}\D_{j}}$ is spectral if and only if $N_j$ divides $b_j$ for each $2\le j\le n$.
\end{theo}

Then $``{\rm(i)}\Rightarrow {\rm(iii)}"$ follows from Theorem \ref{th1.3} and Theorem \ref{th1.4}.
We observe that if $b_n=r_nN_n$ for some integer $r_n$, then
$$\D_{n}\oplus N_n\{0, 1, \ldots, r_n-1\}=\{0, 1, \ldots, b_n-1\}.$$
Here the direct sum $A\oplus B$ means that $a+b$ are all distinct elements for all $a\in{A}$ and $b\in B$. This can yield $``{\rm(iii)}\Rightarrow {\rm(ii)}"$ in Theorem \ref{th1.2}.


A digit set $\D$ is called an {\it integer tile} if $\D$ tiles some cyclic groups $\Z_n$. i.e. there exists $\B$ such that $\D\oplus\B\equiv \Z_n$ (mod $n$). The study of the integer tiles has a long history related to
the geometry of numbers (\cite{CM1999,LL2023,N1997,T1995,S1979} and the references therein). In following Theorem \ref{th1.2}, we present an intriguing result regarding the relationship between the Cantor-Moran spectral measure and integer tile. Write
$$
\mu_{n}=\delta_{\frac{1}{b_1}\D_{1}}\ast\delta_{\frac{1}{b_1b_2}\D_{2}}\ast\cdots\ast\delta_{\frac{1}{b_1b_2\cdots b_n}\D_{n}}=\delta_{\frac{1}{b_1b_2\cdots b_n}{\bf D}_n},
$$
where ${\bf D}_n=\D_{n}+b_n\D_{n-1}+\cdots+b_2\cdots b_n\D_{1}$ is the first $n$ terms iterated digit set of $\{\D_{n}\}_{n=1}^{\infty}$.


\begin{theo}\label{th1.5}
Suppose ${\bf b}=\{b_n\}_{n=1}^{\infty}$ is a sequence of integers bigger than 1 and  ${\bf D}=\{\D_{n}\}_{n=1}^{\infty}$ is a sequence of consecutive digit sets. Then Cantor-Moran measure $\mu_{{\bf b},{\bf D}}$ is spectral if and only if for each $n\geq1$, ${\bf D}_n=\D_{n}\oplus b_n\D_{n-1}\oplus\cdots\oplus b_2\cdots b_n\D_{1}$ is an integer tile.
\end{theo}

We organize this paper as follows. In Section 2, we will study the weak convergence of infinite convolutions and give the proofs of
Theorem \ref{th1.1}. Also, we will introduce some basic definitions, properties of spectral measures. In Section 3, we will discuss the distribution of any bi-zero set of spectral measure $\mu_{\textbf{B},\textbf{D}}$, and   Theorem \ref{th1.3} will be proved. We will devote Section 4 to prove Theorem \ref{th1.2} and Theorem \ref{th1.4}. In Section 5, we prove Theorem \ref{th1.5} and propose an open question on relationship between the spectral Cantor-Moran measure and the tiling of integers.

\section{Notations and Preliminaries}

\subsection{Weak convergence of convolutions}

\

Using Kolmogorov's three series theorem,  Li, Miao and Wang \cite{LMW2022} gave a sufficient and necessary conditions for the existence of infinite convolutions.

\begin{theo}\label{theo3.1}
  Let $\{A_n\}_{n=1}^{\infty}$ be a sequence of nonnegative finite subsets of $\R$ satisfying that $\# A_n\ge2$ for each $n\ge1$. Let $\nu_n=\delta_{A_1}\ast\cdots\ast\delta_{A_n}$. Then the sequence of convolutions $\{\nu_n\}$ converges weakly to a Borel probability measure if and only if
  \begin{equation}\label{eq-Add1}
    \sum_{n=1}^{\infty}\frac{1}{\# A_n}\sum_{a\in A _n}\frac{a}{1+a}<\infty.
  \end{equation}
\end{theo}

\noindent{\bf{ Proof of Theorem \ref{th1.1}}.}
We know from Theorem \ref{theo3.1} that the weak convergence of $\{\mu_n\}_{n=1}^{\infty}$ is equivalent to
\begin{equation}\label{eq3.2}
  \sum_{n=1}^{\infty}\frac{1}{N_n}\sum_{d\in\D_{n}}\frac{d}{b_1\cdots b_n+d}<\infty.
\end{equation}
Note that
\begin{eqnarray}\label{eq3.3}
  &~& \sum_{n=1}^{\infty}\frac{1}{N_n}\sum_{d\in\D_{n}}\frac{d}{b_1\cdots b_n+d} \nonumber\\
   &=& \sum_{\{n:N_n-1>b_1\cdots b_n\}}\frac{1}{N_n}\sum_{d=0}^{N_n-1}\frac{d}{b_1\cdots b_n+d}
   +\sum_{\{n:N_n-1\le b_1\cdots b_n\}}\frac{1}{N_n}\sum_{d=0}^{N_n-1}\frac{d}{b_1\cdots b_n+d}.
\end{eqnarray}

\bigskip

For the sufficiency. Suppose that \eqref{eq1.1} holds. We first claim that $\{n:N_n-1>b_1\cdots b_n\}$ is a finite set. Indeed, if $\{n:N_n-1>b_1\cdots b_n\}$ is an infinite set, then
$$\sum_{n=1}^{\infty}\frac{N_n}{b_1\cdots b_n}\ge\sum_{\{n:N_n-1>b_1\cdots b_n\}}\frac{N_n}{b_1\cdots b_n}=\infty.$$
We get a contradiction. The the claim follows. Hence
\begin{eqnarray}\label{eq-3.5}
\sum_{\{n:N_n-1>b_1\cdots b_n\}}\frac{1}{N_n}\sum_{d=0}^{N_n-1}\frac{d}{b_1\cdots b_n+d}<\infty.
\end{eqnarray}
Notice that
\begin{eqnarray}\label{eq-3.5.1}
&& \sum_{\{n:N_n-1\le b_1\cdots b_n\}}\frac{1}{N_n}\sum_{d=0}^{N_n-1}\frac{d}{b_1\cdots b_n+d}\le\sum_{\{n:N_n-1\le b_1\cdots b_n\}}\frac{1}{N_n}\sum_{d=0}^{N_n-1}\frac{d}{b_1\cdots b_n} \nonumber\\
&=& \sum_{\{n:N_n-1\le b_1\cdots b_n\}}\frac{N_n-1}{2b_1\cdots b_n}<\sum_{n=1}^{\infty}\frac{N_n}{b_1\cdots b_n}<\infty.
\end{eqnarray}
Then \eqref{eq3.2} follows from \eqref{eq3.3}, \eqref{eq-3.5} and \eqref{eq-3.5.1}.

\bigskip
Next we will prove the necessity. Suppose that \eqref{eq3.2} holds. Then it follows from \eqref{eq3.3} that
\begin{equation}\label{eq3.6}
  \sum_{\{n:N_n-1>b_1\cdots b_n\}}\frac{1}{N_n}\sum_{d=0}^{N_n-1}\frac{d}{b_1\cdots b_n+d}<\infty,\,\,\,\sum_{\{n:N_n-1\le b_1\cdots b_n\}}\frac{1}{N_n}\sum_{d=0}^{N_n-1}\frac{d}{b_1\cdots b_n+d}<\infty.
\end{equation}
We assert that $\{n:N_n-1>b_1\cdots b_n\}$ is a finite set. Otherwise,
\begin{eqnarray*}
  \infty &>& 4\sum_{n=1}^{\infty}\frac{1}{N_n}\sum_{d=0}^{N_n-1}\frac{d}{b_1\cdots b_n+d} \\
   &\ge& 4\sum_{\{n:N_n-1>b_1\cdots b_n\}}\frac{1}{N_n}\frac{1}{b_1\cdots b_n+N_n-1}\sum_{d=0}^{N_n-1}d \\
   &=& 4\sum_{\{n:N_n-1>b_1\cdots b_n\}}\frac{2(N_n-1)}{b_1\cdots b_n+N_n-1}=\infty.
\end{eqnarray*}
This make a contradiction, and the assertion follows. Consequently,
\begin{equation}\label{eq-L1}
  \sum_{\{n:N_n-1>b_1\cdots b_n\}}\frac{N_n}{b_1\cdots b_n}<\infty.
\end{equation}
Note that
\begin{eqnarray*}
   && \sum_{\{n:N_n-1\le b_1\cdots b_n\}}\frac{1}{N_n}\sum_{d=0}^{N_n-1}\frac{d}{b_1\cdots b_n+d}\ge\sum_{\{n:N_n-1\le b_1\cdots b_n\}}\frac{1}{N_n}\frac{1}{b_1\cdots b_n+N_n-1}\sum_{d=0}^{N_n-1}d  \nonumber\\
   &=& \sum_{\{n:N_n-1\le b_1\cdots b_n\}}\frac{N_n-1}{2(b_1\cdots b_n+N_n-1)}
   \ge\frac{1}{2}\sum_{\{n:N_n-1\le b_1\cdots b_n\}}\frac{N_n-1}{2b_1\cdots b_n}.
\end{eqnarray*}
It follows that
\begin{equation}\label{eq-L2}
  \sum_{\{n:N_n-1\le b_1\cdots b_n\}}\frac{N_n}{b_1\cdots b_n}<\infty.
\end{equation}
This together with \eqref{eq-L1} yields that
$$\sum_{n=1}^{\infty}\frac{N_n}{b_1\cdots b_n}=\sum_{\{n:N_n-1>b_1\cdots b_n\}}\frac{N_n}{b_1\cdots b_n}+\sum_{\{n:N_n-1\le b_1\cdots b_n\}}\frac{N_n}{b_1\cdots b_n}<\infty.$$
Now we complete the proof.
\qed

\bigskip

\begin{coro}\label{coro1.1}
Let ${\bf b}=\{b_n\}_{n=1}^{\infty}$ be a sequence of integers bigger than 1 and  ${\bf D}=\{\D_{n}\}_{n=1}^{\infty}$ be a sequence of consecutive digit sets with $\D_{n}=\{0,1,\ldots, N_n-1\}$ where $N_n\geq2$. If $N_n\le b_n$ for each $n\geq2$,
then $\mu_n$ converges weakly to a Borel probability measure $\mu_{{\bf b}, {\bf D}}$.
\end{coro}

\begin{proof} As $2\le N_n\le b_n$ for each $n\geq2$, we have
\begin{eqnarray*}
  \sum_{n=1}^{\infty}\frac{N_n}{b_1\cdots b_n}&=& \frac{N_1}{b_1}+\sum_{n=2}^{\infty}\frac{N_n}{b_1\cdots b_n} \\
   &\le& \frac{N_1}{b_1}+\sum_{n=2}^{\infty}\frac{1}{b_1\cdots b_{n-1}} \\
   &\le& \frac{N_1}{b_1}+\sum_{n=2}^{\infty}\frac{1}{2^{n-1}}\\
   &=& \frac{N_1}{b_1}+ 1<\infty.
\end{eqnarray*}
Applying Theorem \ref{th1.1}, the assertion follows.
\end{proof}

\subsection{Spectral measure theoretic preliminaries}

\

Let $\mu$ be a Borel probability measure with compact support on $\R$. The Fourier transform of $\mu$ is defined as usual,
$$
   \widehat{\mu}(\xi)=\int e^{-2\pi i\xi x}\dmu(x)
$$
for any $\xi\in\R$. We will denote by $\ZZ{\left(\widehat{\mu}\right)}=\{\xi \in \R:\widehat{\mu}(\xi)=0\}$ the zero set of $\widehat{\mu}$. In what follows, $e_{\lambda}$ stands for the exponential function $e^{-2\pi i\lambda x}$. Then for a discrete set $\Lambda\subset\R$, $E(\Lambda)=\{e_{\lambda}:\lambda\in\Lambda\}$ is an orthogonal set of $L^2(\mu)$ if and only if  $\widehat{\mu}(\lambda-\lambda')=0$ for $\lambda\neq\lambda'\in\Lambda$, which is equivalent to
\begin{equation}\label{eq2.1}
    (\Lambda-\Lambda)\setminus\{0\}\subset\ZZ{\left(\widehat{\mu}\right)}.
\end{equation}
In this case, we say that $\Lambda$ is a {\it bi-zero set} of $\mu$. Moreover, $\Lambda$ is called a {\it maximal  bi-zero set}  if it is maximal in $\ZZ(\widehat{\mu})$ to have the set difference property. Since bi-zero sets (or spectra) are invariant under translation, without loss of generality, we always assume that $0\in\Lambda$ in this paper.
For $\xi \in \R$, write
$$
   Q_\La(\xi)=\sum_{\la\in\La}\left|\widehat{\mu}(\xi+\la)\right|^2.
$$
The following criterion is a universal test to decide whether a countable set $\Lambda\subset\R$ is a bi-zero set (a spectrum) of $\mu$ or not.

\begin{theo} [\cite{JP1998, LMW20222}]\label{theo2.1}
  Let $\mu$ be a Borel probability measure, and let $\Lambda\subset \R$ be a countable set. Then

 \indent$(\textup{i})$~~ $\La$ is a bi-zero set of $\mu$ if and only if $Q_\Lambda(\xi)\le 1$ for $\xi\in\R$;

 \indent$(\textup{ii})$~~ $\La$ is a spectrum of $\mu$ if and only if $Q_\Lambda(\xi)\equiv 1$ for $\xi\in\R$;

  \indent$(\textup{iii})$~~ $Q_\Lambda(\xi)$ has an entire analytic extension to $\mathbb{C}$ if $\La$ is a bi-zero set of $\mu$.
\end{theo}

As a simple consequence of Theorem \ref{theo2.1}, the following useful theorem was proved in \cite{DHL2014} and will be used to prove our main result.
\begin{theo}\label{theo2.2}
  Let $\mu=\nu\ast \omega$ be the convolution of two probability measures $\nu$ and $\omega$,  and they are not
  Dirac measures. Suppose that $\Lambda$ is a bi-zero set of $\nu$. Then $\Lambda$ is also a bi-zero set of $\mu$, but it
  cannot be a spectrum of $\mu$.
\end{theo}

\section{Proof of Theorem \ref{th1.3} }

In this section, we intend to complete the proof of Theorem \ref{th1.3}. For simplicity, we denote  ${\bf b}_{m,n}=b_m\cdots b_n$ if $m\leq n$. In particularly, if $m=1$, we just denote ${\bf b}_{m,n}$ by ${\bf b}_{n}$. Recall that $\mu_{{\bf b}, {\bf D}}$ is of form convolution, i.e. ${\mu_{{\bf b}, {\bf D}}}=\mu_n\ast\mu_{>n}$ for any integer $n\geq1$, where
\begin{equation}\label{mu_k}
\mu_n:=\ast_{k=1}^{n}\de_{\frac{1}{{\bf b}_{k}}\D_{k}}, \qquad \mu_{>n}:=\ast_{k=n+1}^{\infty}\de_{\frac{1}{{\bf b}_{k}}\D_{k}}.
\end{equation}
So for any bi-zero set $0\in\Lambda$ of $\mu$, we must have
$$\Lambda\setminus\{0\}\subset(\Lambda-\Lambda)\setminus\{0\}\subset{\mathcal Z}(\widehat{\mu}_{{\bf b}, {\bf D}})={\mathcal Z}(\widehat{\mu}_n)\cup{\mathcal Z}(\widehat{\mu}_{>n}).$$
 We now introduce the definition of {\it suitable decomposition}, which will be frequently used in the following of the paper.
\begin{defi}\label{defi1}
Let $\mu=\nu\ast \omega$ be the convolution of two probability measures $\nu$ and $\omega$. Suppose $0\in\Lambda$ is a spectrum of ${\mu}$ and  $0\in\A\subset \Lambda$ is a maximal bi-zero set of $\nu$.
For any $\alpha\in\A$, let
$$\Lambda_{\alpha}:=\Big\{\lambda\in\Lambda: \lambda-\alpha\in{\mathcal Z}(\widehat{\omega})\setminus{\mathcal Z}(\widehat{\nu}) \Big\}\cup\{\alpha\}.$$
We call the following union
\begin{equation*}\label{decomp}
\Lambda=\bigcup_{\alpha\in \A}\Lambda_\alpha
\end{equation*}
 a {\it suitable decomposition} with respect to $(\nu, \omega)$.
\end{defi}
\textbf{Remark:} In the definition, the maximality of $0\in\A\subset \Lambda$  means that for any $\lambda\in\Lambda\setminus\A$, there is a $\alpha\in\A$ such that $\lambda-\alpha\in {\mathcal Z}(\widehat{\omega})\setminus{\mathcal Z}(\widehat{\nu})$.

\bigskip

The suitable decomposition was first introduced by the first author and Wang  \cite{AW2021}. From the Jorgensen-Pedersen lemma \cite{JP1998}, $\Lambda$ is a spectrum for a probability measure $\mu_{{\bf b}, {\bf D}}$ if and only if
$$\sum_{\lambda\in\Lambda}\left|\widehat{\mu}_{{\bf b}, {\bf D}}(\xi+\lambda)\right|^2=1,\quad \forall\ \xi\in\R.$$
Substituting the suitable decomposition $\Lambda=\bigcup_{\alpha\in \A}\Lambda_\alpha$ with respect to $(\mu_n, \mu_{>n})$ into the above equation, we have
$$1\le \sum_{a\in\A}\sum_{\lambda\in\Lambda}\left|\widehat{\mu}(\xi+\lambda)\right|^2
=\sum_{a\in\A}\sum_{\lambda\in\Lambda_a}\left|\widehat{\mu}_n(\xi+\lambda)\right|^2\left|\widehat{\mu}_{>n}(\xi+\lambda)\right|^2.$$
In this paper, we are going to find a suitable decomposition $\Lambda=\bigcup_{\alpha\in \A}\Lambda_\alpha$ with respect to $(\mu_n, \mu_{>n})$ such that $\A$ is a spectrum of $\mu_n$ and $\Lambda_a$ is a spectrum of $\mu_{>n}$ for each $a\in\A$.

With a direct calculation, we have $\mathcal{Z}(\widehat{\delta}_{\D_k})=\frac{1}{N_k}(\Z\setminus N_k\Z)$. So
\begin{equation}\label{muk}
  \mathcal{Z}(\widehat{\mu}_{n})=\bigcup_{k=1}^{n}\frac{{\bf b}_{k}}{N_k}(\Z\setminus N_k\Z),\quad \mathcal{Z}(\widehat{\mu}_{>n})=\bigcup_{k=n +1}^{\infty}\frac{{\bf b}_{k}}{N_k}(\Z\setminus N_k\Z),
\end{equation}
and
\begin{equation}\label{eq2}
  {\mathcal{Z}(\widehat{\mu}_{{\bf b}, {\bf D}})}=\mathcal{Z}(\widehat{\mu}_{n})\cup\mathcal{Z}(\widehat{\mu}_{>n})
  =\bigcup_{k=1}^{\infty}\frac{{\bf b}_{k}}{N_k}(\Z\setminus N_k\Z).
\end{equation}
Now we give an important observation and it will be used to guarantee that there is a strong link between the spectral measures $\mu_n, \mu_{>n}$ and $\mu_{{\bf b}, {\bf D}}$.

\begin{lem}\label{lem1}
Suppose that $\mu_{{\bf b}, {\bf D}}$ is a spectral measure and $\{\lambda, \gamma\}$ is a bi-zero set of $\mu_{{\bf b}, {\bf D}}$. If $\lambda\in \ZZ(\widehat{\delta}_{\frac{1}{{\bf b}_{n}}\D_{n}})$ and $\gamma\in\ZZ(\widehat{\delta}_{\frac{1}{{\bf b}_{k}}\D_{k}})\setminus \ZZ(\widehat{\delta}_{\frac{1}{{\bf b}_{n}}\D_{n}})$ with $k>n$, then we must have
$$
\lambda-\gamma\in \ZZ(\widehat{\delta}_{\frac{1}{{\bf b}_{n}}\D_{n}}).
$$
\end{lem}
\begin{proof} The bi-zero property of $\{\lambda, \gamma\}$ implies that
$$\lambda-\gamma\in\ZZ(\widehat{\mu}_{{\bf b}, {\bf D}}).$$
Suppose to the contrary that
$$\lambda-\gamma\in\ZZ(\widehat{\delta}_{\frac{1}{{\bf b}_{j}}\D_{j}})\setminus \ZZ(\widehat{\delta}_{\frac{1}{{\bf b}_{n}}\D_{n}})$$
for some $j\neq n$. Then from \eqref{muk}, we can write them as
$$
\lambda=\frac{{\bf b}_{n}}{N_{n}}a_n,\
\gamma=\frac{{\bf b}_{k}}{N_{k}}a_k\ \text{ and }\ \lambda-\gamma=\frac{{\bf b}_{j}}{N_{j}}a_j,
$$
where $a_i\in \Z\setminus N_{i}\Z, i\in\{n, k, j\}$.  Without loss of generality, we assume $j>n$. After some rearrangement, we have
\begin{equation}\label{eq3}
  \frac{a_n}{N_{n}}=\frac{{\bf b}_{n+1,k}}{N_{k}}a_k+\frac{{\bf b}_{n+1,j}}{N_{j}}a_j.
\end{equation}
Reduce all fractions in above equation to be their simplest form, i.e.,
$$\frac{a_n'}{N_{n}'}=\frac{{\bf b}_{n+1,k}'}{N_{k}'}a_k+\frac{{\bf b}_{n+1,j}'}{N_{j}'}a_j=
\frac{{\bf b}_{n+1,k}'a_kN_j'+{\bf b}_{n+1,j}'a_jN_k'}{N_k'N_j'},$$
where
\begin{equation}\label{eq3'}
\gcd(a_n', N_{n}')=1, \gcd({\bf b}_{n+1,k}', N_{k}')=1\text{ and } \gcd({\bf b}_{n+1,j}', N_{j}')=1.
\end{equation}
It implies $N_{n}'$ divides $N_{k}'N_j'$.  Since $a_n\in \Z\setminus N_{n}\Z$, we have $N_{n}'>1$. Let $s_{n}$ be a prime factor of $N_{n}'$ and of course it is also a prime factor of $N_{n}$. Then we must have
$$s_{n}|N_{k}'\text{ or } s_{n}|N_{j}'.$$
Without loss of generality, we assume $N_{k}'=s_{n}t_{k}'$ for some integer $t_k'$. It follows from the second equation in \eqref{eq3'} that
\begin{equation}\label{eq3''}
\gcd(s_{n}, {\bf b}_{n+1,k}')=1.
\end{equation}

\medskip

Denote $t_{k}=\gcd(N_{k}, {\bf b}_{n+1, k})$. Then  $N_{k}=N_k't_k=s_{n}t_{k}t_{k}'$ and ${\bf b}_{n+1, k}={\bf b}_{n+1, k}'t_k$. Let $\mathcal{E}_{t_{k}}=\{0,1,\ldots,t_k-1\}$, $\mathcal{E}_{t_{k}'}=\{0,1,\ldots,t_k'-1\}$ and $\mathcal{E}_{s_n}=\{0,1,\ldots,s_n-1\}$. We can factorize $\D_{k}$ as
$$\D_{k}=\mathcal{E}_{t_{k}}\oplus t_{k}\mathcal{E}_{s_{n}}\oplus s_{n}t_{k}\mathcal{E}_{t_{k}'}.$$
Denote
$$\nu=\ast_{i\neq {k}}\delta_{\frac{1}{{\bf b}_{i}}\D_{i}}\ast \delta_{\frac{1}{{\bf b}_{k}}\mathcal{E}_{t_{k}}}\ast\delta_{\frac{1}{{\bf b}_{k}}s_{n}t_{k}\mathcal{E}_{t_{k}'}}. $$
Then
${\mu_{{\bf b}, {\bf D}}=\nu\ast\delta_{\frac{1}{{\bf b}_{k}}t_{k}\mathcal{E}_{s_{n}}}}.$
Note that
\begin{eqnarray*}
  \mathcal{Z}(\widehat{\delta}_{\frac{1}{{\bf b}_{k}}t_{k}\mathcal{E}_{s_{n}}}) &=& \frac{{\bf b}_{n}{\bf b}_{n+1,k}}{t_{k}s_{n}}(\Z\setminus s_{n}\Z)\\
   &=& \frac{{\bf b}_{n}{\bf b}_{n+1,k}'}{s_{n}}(\Z\setminus s_{n}\Z)\\
  &\subset&\frac{{\bf b}_{n}}{s_{n}}(\Z\setminus s_{n}\Z)\quad \left(\text{ as } \gcd({\bf b}_{n+1,k}', s_{n})=1\right)\\
  &\subset&\frac{{\bf b}_{n}}{N_{n}}(\Z\setminus N_{n}\Z) \quad \left(\text{ as } s_{n}|N_{n}\right)\\
  &=&\mathcal{Z}(\widehat{\delta}_{\frac{1}{{\bf b}_{n}}\D_{n}}) \subset {\ZZ(\widehat{\nu})}.
\end{eqnarray*}
So
$$\ZZ(\widehat{\mu}_{{\bf b}, {\bf D}})=\ZZ(\widehat{\nu})\cup \mathcal{Z}(\widehat{\delta}_{\frac{1}{{\bf b}_{k}}t_{k}\mathcal{E}_{s_{n}}})=\ZZ(\widehat{\nu}).$$
Let $\Lambda$ be an arbitrary bi-zero set of $\mu_{{\bf b}, {\bf D}}$, the above equation implies that it also a bi-zero set of $\nu$. It follows from Theorem \ref{theo2.2} that $\Lambda$ can not be a spectrum of $\mu_{{\bf b}, {\bf D}}$. Therefore, $\mu_{{\bf b}, {\bf D}}$ is not a spectral measure, a contradiction. The proof of the lemma is now completed.
\end{proof}

\begin{lem}\label{lem2}
Suppose that {$\mu_{{\bf b}, {\bf D}}$} is a spectral measure and $\{\lambda, \gamma\}$ is a bi-zero set of $\mu_{{\bf b}, {\bf D}}$. Let $n\geq1$ be an integer. \\
\rm{(i)} If $\lambda\in \ZZ(\widehat{\mu}_n)$ and $\gamma\in {\ZZ(\widehat{\mu}_{>n})\setminus \ZZ(\widehat{\mu}_n)}$, then we must have
$$
\lambda-\gamma\in \ZZ(\widehat{\mu}_n).
$$
\rm{(ii)} If $\lambda, \gamma\in {\ZZ(\widehat{\mu}_{>n})\setminus \ZZ(\widehat{\mu}_n)}$, then we must have
$$
\lambda-\gamma\in \ZZ(\widehat{\mu}_{>n})\setminus \ZZ(\widehat{\mu}_n).
$$
\end{lem}
\begin{proof}
\rm{(i)} The assumption and $\eqref{muk}$ imply that
$\lambda\in\ZZ(\widehat{\delta}_{\frac{1}{{\bf b}_{k}}\D_{k}})$ and $\gamma\in\ZZ(\widehat{\delta}_{\frac{1}{{\bf b}_{m}}\D_{m}})\setminus\ZZ(\widehat{\delta}_{\frac{1}{{\bf b}_{k}}\D_{k}})$
for some $k\le n<m$. From Lemma \ref{lem1} we have
$$\lambda-\gamma\in\ZZ(\widehat{\delta}_{\frac{1}{b_{k}}\D_{k}})\subset \ZZ(\widehat{\mu}_n).$$

\rm{(ii)} It should be noticed that $\ZZ(\widehat{\nu})=-\ZZ(\widehat{\nu})$ for any  measure $\nu$.
Suppose to the contrary that
$$\lambda':=\lambda-\gamma\in\ZZ(\widehat{\mu}_n).$$
Since $$\lambda=\lambda'-(-\gamma)\in\mathcal{Z}(\widehat{\mu}_{>n})\setminus \ZZ(\widehat{\mu}_n)\subset\ZZ(\widehat{\mu}_{{\bf b}, {\bf D}}),$$
$\{\lambda',-\gamma\}$ is a bi-zero set of $\mu_{{\bf b}, {\bf D}}$ with $\la'\in \ZZ(\widehat{\mu}_n)$ and $-\gamma\in \mathcal{Z}(\widehat{\mu}_{>n})\setminus \ZZ(\widehat{\mu}_n)$. It follows from \rm{(i)}  that $\lambda=\lambda'-(-\gamma)\in\ZZ(\widehat{\mu}_n)$, a contradiction. We have completed the proof.
\end{proof}

\bigskip

 As a consequence of Lemma \ref{lem2}, we have

\begin{coro}\label{coro-lem12}
Suppose $0\in\Lambda$ is a spectrum of $\mu_{{\bf b}, {\bf D}}$, then $\Lambda\cap \ZZ(\widehat{\mu}_n)\neq \emptyset$ for any $n\geq1$.
\end{coro}
\begin{proof} The bi-zero property of $\Lambda$ implies that
$$\Lambda\setminus\{0\}\subset(\Lambda-\Lambda)\setminus\{0\}\subset{\mathcal Z}(\widehat{\mu}_{{\bf b}, {\bf D}})={\mathcal Z}(\widehat{\mu}_{n})\cup{\mathcal Z}(\widehat{\mu}_{>n}).$$
If the assertion is not true for some $n\geq1$, then we have
$$\Lambda\setminus\{0\}\subset{\mathcal Z}(\widehat{\mu}_{>n})\setminus {\mathcal Z}(\widehat{\mu}_{n}).$$
From Lemma \ref{lem2} (ii), we have
$$(\Lambda-\Lambda)\setminus\{0\}\subset{\mathcal Z}(\widehat{\mu}_{>n})\setminus {\mathcal Z}(\widehat{\mu}_{n}).$$
That is to say $\Lambda$ is a bi-zero set of $\mu_{>n}$. It follows from Theorem \ref{theo2.2} that $\Lambda$ cannot be a spectrum of $\mu_{{\bf b}, {\bf D}}$, which is a contradiction. We have completed the proof.
\end{proof}

\bigskip

\begin{lem}\label{lem3}
 Let $0\in\Lambda$ be a spectrum of ${\mu}_{{\bf b}, {\bf D}}$, then $\Lambda\subset \frac{1}{N_1}\Z$.
\end{lem}
\begin{proof} We write ${\mu_{{\bf b}, {\bf D}}}=\mu_1\ast\mu_{>1}$ and decompose
$
\Lambda=\bigcup_{\alpha\in \A}\Lambda_\alpha
$
into a suitable form with respect to $(\mu_1, \mu_{>1})$.  That is to say $0\in\A\subset \Lambda$ is a maximal bi-zero set of ${\mu}_{1}$ and $$\Lambda_{\alpha}:=\Big\{\lambda\in\Lambda: \lambda-\alpha\in{\mathcal Z}(\widehat{\mu}_{>1})\setminus{\mathcal Z}(\widehat{\mu}_{1})\Big\}\cup\{\alpha\}, \forall \alpha\in\A.$$
It follows  from Corollary \ref{coro-lem12} that $\A\setminus\{0\}$ is not empty. Recall that
\begin{eqnarray*}
\ZZ(\widehat{\mu}_{1})&=&\ZZ(\widehat{\delta}_{\frac{1}{{\bf b}_1}\D_{1}})
=\frac{{\bf b}_1}{N_1}(\Z\setminus N_1\Z)
\subset \frac{1}{N_1}\Z.
\end{eqnarray*}
For any element $\lambda\in\Lambda$, if $\lambda\in\A$, then
$$\lambda\in\A\subset\ZZ(\widehat{\mu}_{1})\cup\{0\}\subset \frac{1}{N_1}\Z.$$
Otherwise, the maximality of $\A$ implies that there is a $\alpha\in\A$ such that
$$\lambda-\alpha\in{\mathcal Z}(\widehat{\mu}_{>1})\setminus{\mathcal Z}(\widehat{\mu}_{1}).$$
Take an element $\alpha'\in\A\setminus\{\alpha\}$, then $\alpha'-\alpha\in{\mathcal Z}(\widehat{\mu}_{1}).$ From Lemma \ref{lem2} (i), we have
$$\alpha'-\lambda=(\alpha'-\alpha)-(\lambda-\alpha)\in{\mathcal Z}(\widehat{\mu}_{1}).$$
Hence
$$\lambda=\alpha'-(\alpha'-\lambda)\in \A-{\mathcal Z}(\widehat{\mu}_{1})\subset \frac{1}{N_1}\Z.$$
So $\Lambda\subset \frac{1}{N_1}\Z$.
\end{proof}

\bigskip

 Theorem \ref{th1.3} shows that if ${\mu}_{{\bf b}, {\bf D}}$ is a spectral measure, then its any ``truncation'' is still a spectral measure. The proof is inspired by \cite[Proposition 4.3]{AW2021}.

\noindent{\bf{ Proof of Theorem \ref{th1.3}}.}
We first claim that the $\{\Lambda_{\alpha}\}_{\alpha\in\A}$ are disjoint pairwise. Otherwise, suppose there is an element $\lambda$ such that
$$\lambda\in\Lambda_{\alpha}\cap\Lambda_{\alpha'}$$
for some $\alpha\neq\alpha'\in\A$. Since $\alpha-\alpha'\in{\mathcal Z}(\widehat{\mu}_{n})$, we have $\alpha\in\Lambda_{\alpha}\setminus\Lambda_{\alpha'}$ and $\alpha'\in\Lambda_{\alpha'}\setminus\Lambda_{\alpha}$. So $\lambda\neq\alpha$ and $\lambda\neq\alpha'$.
From the definition of $\Lambda_{\alpha}$, we have
$$\lambda-\alpha, \lambda-\alpha'\in{\mathcal Z}(\widehat{\mu}_{>n})\setminus{\mathcal Z}(\widehat{\mu}_{n}).$$
Taking $\lambda_1=\lambda-\alpha', \lambda_2=\lambda-\alpha$ in Lemma \ref{lem2}, we have
$$\alpha-\alpha'=(\lambda-\alpha')-(\lambda-\alpha)\in{\mathcal Z}(\widehat{\mu}_{>n})\setminus{\mathcal Z}(\widehat{\mu}_{n}).$$
It contradicts to the fact that $\A$ is a bi-zero set of $\mu_n$. The claim holds.

\bigskip
For any $\alpha\in\A$ and  $\lambda_{\alpha}\neq\lambda_{\alpha}'\in\Lambda_\alpha$, we have
$$\lambda_{\alpha}-\alpha,\ \lambda_{\alpha}'-\alpha\in {\mathcal Z}(\widehat{\mu}_{>n})\setminus{\mathcal Z}(\widehat{\mu}_{n}).$$
From  Lemma \ref{lem2}, we have
\begin{equation}\label{eq-alpha}
\lambda_\alpha-\lambda_\alpha'=(\lambda_\alpha-\alpha)-(\lambda_\alpha'-\alpha)
\in{\mathcal Z}(\widehat{\mu}_{>n})\setminus{\mathcal Z}(\widehat{\mu}_{n}).
\end{equation}
For any $\alpha'\neq\alpha\in\A$, as $\lambda_{\alpha}\not\in\Lambda_{\alpha'}$, we have
$$\lambda_\alpha-\alpha'\in{\mathcal Z}(\widehat{\mu}_{n}).$$
As $\lambda_{\alpha'}-\alpha'\in{\mathcal Z}(\widehat{\mu}_{>n})\setminus{\mathcal Z}(\widehat{\mu}_{n})$ for any $\lambda_{\alpha'}\in\Lambda_{\alpha'}$, then Lemma \ref{lem1} implies that
\begin{equation}\label{eq-alpha'}
\lambda_\alpha-\lambda_{\alpha'}=(\lambda_\alpha-\alpha')-(\lambda_{\alpha'}-\alpha')\in{\mathcal Z}(\widehat{\mu}_{n}).
\end{equation}
Summarizing equation \eqref{eq-alpha} and \eqref{eq-alpha'}, we have
\begin{equation}\label{eq.add.4.1.0}
   (\Lambda_\alpha-\Lambda_\alpha)\setminus\{0\}\subset{\mathcal Z}(\widehat{\mu}_{>n})\setminus{\mathcal Z}(\widehat{\mu}_{n}) \quad \text{and} \quad  \Lambda_\alpha-\Lambda_{\alpha'}\subset{\mathcal Z}(\widehat{\mu}_{n}).
\end{equation}


For any $\xi\in[0,1]$, we have
\begin{equation}\label{eq.add.4.1.1}
Q_{\Lambda}(\xi)=\sum_{\lambda\in\Lambda}|\widehat{\mu}_{{\bf b}, {\bf D}}(\xi+\lambda)
|^2=\sum_{\alpha\in\A}\sum_{\lambda_\alpha\in\Lambda_\alpha}
|\widehat{\mu}_{n}(\xi+\lambda_\alpha)|^2
|\widehat{\mu}_{>n}(\xi+\lambda_\alpha)|^2.
\end{equation}
 We define an equivalence relationship $\sim$  such that $\lambda\sim\lambda'$ whenever $\lambda'-\lambda\in b_1b_2\cdots b_n\Z$.
Then the quotient group
$\Lambda_{\alpha}/\sim=\{[\lambda]: \lambda\in\Lambda_{\alpha}\}$ is a partition of $\Lambda_\alpha$, where
$$[\lambda]:=\{\lambda'\in\Lambda_{\alpha}: \lambda'\sim\lambda\}.$$
Since $\Lambda\subset \frac{1}{N_1}\Z$ (from Lemma \ref{lem3}), $\Lambda_{\alpha}/\sim$ is a finite set, denote it by $\{[\lambda_{\alpha, 1}], \ldots, [\lambda_{\alpha, n_{\alpha}}]\}$.
Note that  $\widehat{\mu}_{n}$ is $b_1b_2\cdots b_n$-periodic. For any $\lambda\in[\lambda_{\alpha, i}]$, we have
$$|\widehat{\mu}_{n}(\xi+\lambda)|=|\widehat{\mu}_{n}(\xi+\lambda_{\alpha,i})|,\quad \forall\ \xi\in\R .$$
For any $\xi\in[0, 1]$ and $\alpha\in\A$, there is the unique $\lambda_{\alpha, i(\xi)},i(\xi)\in\{1, 2, \ldots, n_{\alpha}\},$ such that
\begin{eqnarray*}
|\widehat{\mu}_{n}(\xi+\lambda_{\alpha, i(\xi)})|&=&\max\Big\{|\widehat{\mu}_{n}(\xi+\lambda_{\alpha, i})|:1\le i\le n_{\alpha}\Big\}\\
&=&\max\Big\{|\widehat{\mu}_{n}(\xi+\lambda_\alpha)|:
\lambda_\alpha\in\Lambda_\alpha\Big\}.
\end{eqnarray*}
As $\A$ and $\Lambda_\alpha/\sim$ are finite set, we can find  a finite set $\{\lambda_{\alpha, i_{\alpha}}\}_{\alpha\in\A}$ such that $\lambda_{\alpha, i(\xi_j)}=\lambda_{\alpha,i_\alpha}$ for infinitely many $\{\xi_{j}\}_{j=1}^{\infty}$. Combined with \eqref{eq.add.4.1.1}, we have
\begin{eqnarray*}\label{eq.add.4.1.2}
Q_{\Lambda}(\xi_{j})&\leq&\sum_{\alpha\in\A}
|\widehat{\mu}_{n}(\xi_{j}+\lambda_{\alpha,i_\alpha})|^2
\Big(\sum_{\lambda_\alpha\in\Lambda_\alpha}
|\widehat{\mu}_{>n}(\xi_{j}+\lambda_\alpha)|^2\Big)
\nonumber\\
&\leq& \sum_{\alpha\in\A}
|\widehat{\mu}_{n}(\xi_{j}+\lambda_{\alpha,i_\alpha})|^2\nonumber\\
&\leq&1,
\end{eqnarray*}
where the last two inequalities follow from  \eqref{eq.add.4.1.0}  and Theorem \ref{theo2.1} (i).
On the other hand, as $\Lambda$ is a spectrum of ${\mu}_{{\bf b}, {\bf D}}$, we have $Q_{\Lambda}(\xi)\equiv1$. This forces
$$\sum_{\lambda_\alpha\in\Lambda_\alpha}
|\widehat{\mu}_{>n}(\xi_{j}+\lambda_\alpha)|^2\equiv1\quad \text{and} \quad \sum_{\alpha\in\A}
|\widehat{\mu}_{n}(\xi_{j}+\lambda_{\alpha,i_\alpha})|^2\equiv1,\quad\ \forall~j\geq1.$$
As all $\xi_j\in[0, 1]$, the entire function property implies that
$$\sum_{\lambda_\alpha\in\Lambda_\alpha}
|\widehat{\mu}_{>n}(\xi+\lambda_\alpha)|^2\equiv1\quad \text{and} \quad  \sum_{\alpha\in\A}
|\widehat{\mu}_{n}(\xi+\lambda_{\alpha,i_\alpha})|^2\equiv1,\quad\ \forall~\xi\in\R.$$
Hence $\A_0=\{\lambda_{\alpha,i_{\alpha}}\}_{\alpha\in\A}\subset\Lambda$ is a spectrum of ${\mu}_{n}$ and each $\Lambda_{\alpha}$ is a spectrum of $\mu_{>n}$. The desired result follows.
\qed

\section{Proofs of Theorem \ref{th1.2} and Theorem \ref{th1.4}}

Theorem \ref{th1.3} gives a necessary condition for $\mu_{{\bf b}, {\bf D}}$ to be spectral. In this section, we will focus on analysing the spectrality of $\mu_n, n\geq2$. From equation \eqref{muk}, for any $1\le k<n$ we have
$$\mu_n=\mu_{k}\ast\mu_{k+1, n},$$
where $\mu_k=\ast_{j=1}^{k}\delta_{\frac{1}{{\bf b}_j}\D_{j}}$ and $\mu_{k+1,n}=\ast_{j=k+1}^{n}\delta_{\frac{1}{{\bf b}_j}\D_{j}}$.
With the same proof, the conclusions in Lemma \ref{lem2} and Theorem \ref{th1.3} are also true for $\mu_n$.

\begin{lem}\label{coro1}
Suppose that $\mu_{n}$ is a spectral measure and $\{\lambda, \gamma\}$ is a bi-zero set of $\mu_{n}$ where $n\geq2$. Let $1\le k<n$ be an integer. \\
\rm{(i)} If $\la\in \ZZ(\widehat{\mu}_{k})$ and $\gamma\in\ZZ(\widehat{\mu}_{k+1,n})\setminus \ZZ(\widehat{\mu}_{k})$, then we must have
$
\la-\gamma\in \ZZ(\widehat{\mu}_{k}).
$
\\
\rm{(ii)} If $\lambda, \gamma\in\ZZ(\widehat{\mu}_{k+1,n})\setminus \ZZ(\widehat{\mu}_{k})$, then we must have
$\lambda-\gamma\in \ZZ(\widehat{\mu}_{k+1,n})\setminus \ZZ(\widehat{\mu}_{k})$.
\end{lem}
\begin{proof} The proof is just the same as  Lemma \ref{lem2}.
\end{proof}
\bigskip

\begin{theo}\label{coro2}
Suppose $0\in\Lambda$ is a spectrum of $\mu_n$. There is a suitable decomposition $\Lambda=\bigcup_{\alpha\in\A}\Lambda_{\alpha}$ with respect to  $(\mu_{k}, {\mu}_{k+1,n})$ such that $\A$ is a spectrum of $\mu_k$ and $\Lambda_{\alpha}$ is a spectrum of $\mu_{k+1, n}$ for each $a\in\A$.
\end{theo}
\begin{proof}
The proof is just the same as Theorem \ref{th1.3}.
\end{proof}

\medskip

\begin{lem}\label{lemD_q}
Suppose $\mathcal{E}_{N}=\{0, 1, \ldots, N-1\}$ is a consecutive digit set with cardinality $N$. Then $\delta_{\mathcal{E}_N}$ is a spectral measure. Moreover,  $0\in\C$  is a spectrum of $\delta_{\mathcal{E}_N}$ if and only if $\#\C=N$ and $\C\equiv \left\{0, \frac1N, \ldots, \frac{N-1}{N}\right\}\pmod{\Z}$.
\end{lem}
\begin{proof} It has been proved in \cite{AH2014} that $\delta_{\mathcal{E}_N}$ is a spectral measure and admits a spectrum $\left\{0, \frac1N, \cdots, \frac{N-1}{N}\right\}$. Since $\widehat{\delta}_{\mathcal{E}_N}$ is $1-$period, for any digit set $0\in\C$ with $\#\C=N$ and $\C\equiv \left\{0, \frac1N, \cdots, \frac{N-1}{N}\right\}\pmod{\Z}$, we have
$$\sum_{c\in\C}\left|\widehat{\delta}_{\mathcal{E}_N}(\xi+c)\right|^2=\sum_{j=0}^{N-1}\left|
\widehat{\delta}_{\mathcal{E}_N}(\xi+\frac{j}{N})\right|^2.$$
It follows from Theorem \ref{theo2.1} (ii) that $\C$ is a spectrum of ${\delta}_{\mathcal{E}_N}$.

\medskip
Next we will prove the sufficiency. Suppose $0\in\C$ is a spectrum of $\delta_{\mathcal{E}_N}$.
The bi-zero property implies that
$$\C\setminus\{0\}\subset(\C-\C)\setminus\{0\}\subset\ZZ(\widehat{\delta}_{\mathcal{E}_N})
=\frac{1}{N}(\Z\setminus N\Z).$$
So $\C\pmod{\Z}\subset \left\{0, \frac1N, \ldots, \frac{N-1}{N}\right\}.$ The completeness implies $\#\C=\dim L^2(\delta_{\mathcal{E}_N})=N$, so
$$\C\equiv \left\{0, \frac1N, \ldots, \frac{N-1}{N}\right\}\pmod{\Z}.$$
\end{proof}

\medskip

\noindent{\bf{ Proof of Theorem \ref{th1.4}}.} The sufficiency can be found in \cite{AH2014}. For the necessity, suppose on the contrary that there is $2\le k\le n$ such that $N_k$ does not divide $b_k$.  Denote  $d={\gcd(b_k, N_k)}$ and write $b_k=db_k', N_k=dN_k'$. Then $N_k'$ has a prime factor $s_k$ which does not divide $b_k'$. Let $N_k'=s_kt_k$, then $N_k=dt_k s_k$.

 According to Theorem \ref{coro2}, the spectrality of $\mu_{n}=\mu_{k-2}\ast\mu_{k-1,n}$ implies that $\mu_{k-1, n}$ is also a spectral measure. If $n=k$, then $\mu_{k-1, n}=\mu_{k-1,n}$. If $n>k$, then $\mu_{k-1,n}=\mu_{k-1,k}\ast\mu_{k+1,n}$. From Theorem \ref{coro2} again, $\mu_{k-1, k}$ is a spectral measure too.  Suppose $0\in\Lambda$ is a spectrum of $\mu_{k-1, k}=\delta_{\frac{1}{{\bf b}_{k-1}}\D_{k-1}}\ast\delta_{\frac{1}{{\bf b}_k}\D_{k}}$,  we can decompose it into a suitable form
$$\Lambda=\bigcup_{a\in\A}\Lambda_{a}$$
with respect to $(\delta_{\frac{1}{{\bf b}_{k-1}}\D_{k-1}}, \delta_{\frac{1}{{\bf b}_k}\D_{k}})$, that is to say  $0\in \A$ is a maximal set of $\delta_{\frac{1}{{\bf b}_{k-1}}\D_{k-1}}$ and each $$\Lambda_{a}=\left\{\lambda\in\Lambda: \lambda-a\in\ZZ(\widehat{\delta}_{\frac{1}{{\bf b}_{k}}\D_{k}})\setminus
\ZZ(\widehat{\delta}_{\frac{1}{{\bf b}_{k-1}}\D_{k-1}})\right\}\cup\{a\}$$ is a spectrum of $\delta_{\frac{1}{{\bf b}_k}\D_{k}}$. Moreover, one can check that $\{\Lambda_{a}\}_{a\in\A}$ are disjoint pairwise. As $0\in\Lambda_0$, it follows from Lemma \ref{lemD_q} that
 $$\Lambda_0={\bf b}_k\left\{0, \frac{1}{N_k}+m_1, \cdots, \frac{N_k-1}{N_k}+m_{N_k-1}\right\}.$$
Fix an element $a\in\A\setminus\{0\}\subset\ZZ(\widehat{\delta}_{\frac{1}{{\bf b}_{k-1}}\D_{k-1}})$, then $a=\frac{{\bf b}_{k-1}}{N_{k-1}}m$ for some $m\in\Z\setminus N_{k-1}\Z$. Recall that $N_k=dt_ks_k$ where $s_k$ is a prime integer and $1\le t_k\le N_k-1$. So
$${\bf b}_k\left(\frac{t_k}{N_k}+m_{t_k}\right)\in\Lambda_0.$$
It together with $\Lambda_0\cap\Lambda_{a}=\emptyset$ implies that
${\bf b}_k\left(\frac{t_k}{N_k}+m_{t_k}\right)\not\in\Lambda_a$ and so
$$\frac{{\bf b}_{k-1}}{N_{k-1}}m- {\bf b}_k\left(\frac{t_k}{N_k}+m_{t_k}\right)=\frac{{\bf b}_{k-1}}{N_{k-1}}m'\in \ZZ(\widehat{\delta}_{\frac{1}{{\bf b}_{k-1}}\D_{k-1}}),$$
for some $m'\in\Z\setminus N_{k-1}\Z$.
After some rearrangement, we have
$$\frac{m-m'-N_{k-1}b_km_{t_k}}{N_{k-1}}=\frac{b_kt_k}{N_k}=\frac{b_k'}{s_k}.$$
It forces $s_k$ divides $N_{k-1}$ since $\gcd(s_k, b_k')=1$. Let $\mathcal{E}_{d}=\{0,1,\ldots,d-1\}$, $\mathcal{E}_{s_k}=\{0,1,\ldots,s_k-1\}$ and $\mathcal{E}_{t_k}=\{0,1,\ldots,t_k-1\}$.
We factorize $\D_{k}$ as
$$\D_{k}=\mathcal{E}_{d}\oplus d\mathcal{E}_{s_{k}}\oplus ds_{k}\mathcal{E}_{t_{k}}.$$
Denote
$$\nu=\delta_{\frac{1}{{\bf b}_{k-1}}\D_{k-1}}\ast \delta_{\frac{1}{{\bf b}_{k}}\mathcal{E}_{d}}\ast\delta_{\frac{1}{{\bf b}_{k}}s_{k}d\mathcal{E}_{t_{k}}}. $$
Then $\mu_{k-1, k}=\nu\ast\delta_{\frac{1}{{\bf b}_{k}}d\mathcal{E}_{s_{k}}}.$
Note that
\begin{eqnarray*}
  \mathcal{Z}(\widehat{\delta}_{\frac{1}{{\bf b}_{k}}d\mathcal{E}_{s_{k}}}) &=& \frac{{\bf b}_{k-1}b_{k}}{ds_{k}}(\Z\setminus s_{k}\Z)=\frac{{\bf b}_{k-1}b_{k}'}{s_{k}}(\Z\setminus s_{k}\Z)\\
  &\subset&\frac{{\bf b}_{k-1}}{s_{k}}(\Z\setminus s_{k}\Z)\quad \quad\left(\text{ as } \gcd(b_{k}', s_{k})=1\right)\\
  &\subset&\frac{{\bf b}_{k-1}}{N_{k-1}}(\Z\setminus N_{k-1}\Z) \quad \quad  \left(\text{ as } s_{k}|N_{k-1}\right)\\
  &=&\mathcal{Z}(\widehat{\delta}_{\frac{1}{{\bf b}_{k-1}}\D_{k-1}}) \subset\ZZ(\widehat{\nu}).
\end{eqnarray*}
So
$$\ZZ(\widehat{\mu}_{k-1, k})=\ZZ(\widehat{\nu})\cup \mathcal{Z}(\widehat{\delta}_{\frac{1}{{\bf b}_{k}}d\mathcal{E}_{s_{k}}})=\ZZ(\widehat{\nu}).$$
This implies that $\Lambda$ is also a bi-zero set of $\nu$. It follows from Theorem \ref{theo2.2} that $\Lambda$ can not be a spectrum of $\mu_{k-1, k}$, which is a contradiction. Hence $N_k$ divides $b_k$ for each $2\le k\le n$. We have completed the proof.
\qed

Before proving Theorem \ref{th1.2}, we need some useful lemmas. The following lemma is well-known. For the reader's convenience, we give a proof.
\begin{lem}\label{lem4}
 $\Lambda$ is a spectrum of $\mu_{{\bf b},{\bf D}}$ if and only if $\frac{1}{a}\Lambda$ is a spectrum of $\mu_{{\bf b},{\bf D}}$ for any $a\in\R\setminus\{0\}$.
\end{lem}
\begin{proof} Note that
$$\widehat{\mu}_{{\bf b},{\bf D}}(\xi)=\prod_{n=1}^{\infty}\widehat{\delta}_{\frac{1}{{\bf b}_n}a\D_{n}}(\xi)=
\prod_{n=1}^{\infty}\widehat{\delta}_{\frac{1}{{\bf b}_n}\D_{n}}(a\xi)=
\widehat{\mu}_{{\bf b},{\bf D}}(a\xi).$$
Hence
$$\sum_{\gamma\in \frac{1}{a}\Lambda}\left|\widehat{\mu}_{{\bf b},{\bf D}}(\xi+\gamma)\right|^2=\sum_{\gamma\in \frac{1}{a}\Lambda}\left|\widehat{\mu}_{{\bf b},{\bf D}}(a(\xi+\gamma))\right|^2
=\sum_{\lambda\in\Lambda}\left|\widehat{\mu}_{{\bf b},{\bf D}}(a\xi+\lambda)\right|^2.$$
The assertion follows from Theorem \ref{theo2.1}.
\end{proof}

The following theorem has been proved by Gabardo and Lai \cite{GL2014}.

\begin{theo}[\cite{GL2014}]\label{th-Add1}
Any positive Borel measure $\mu$ and $\nu$ such that $\mu\ast\nu=\mathcal{L}_{[0,1]}$ are spectral measures.
\end{theo}

Now we have all ingredients for the proof of Theorem \ref{th1.2}.

\noindent{\bf{ Proof of Theorem \ref{th1.2}}.}
$``{\rm(i)}\Rightarrow {\rm(iii)}"$ Suppose $N_n$ divides $b_n$ for each $n\geq2$. Let $b_1'=N_1$ and $b_n'=b_n$ for $n\geq2$. Then
${\bf b}'=\{b_n'\}_{n=1}^{\infty}$ is a sequence of integers  satisfies $N_n|b_n'$ for each $n\geq1$. From \cite[Theorem 1.4]{AH2014}, we have $\mu_{{\bf b}',{\bf D}}=\mu_{{\bf b},\frac{b_1}{N_1}{\bf D}}$ is a spectral measure. It follows from Lemma \ref{lem4} that $\mu_{{\bf b},{\bf D}}$ is a spectral measure.

$``{\rm(iii)}\Rightarrow {\rm(ii)}"$ Suppose $\mu_{{\bf b},{\bf D}}$ is a spectral measure, according to Theorem \ref{th1.3} and Theorem \ref{th1.4}, we have $N_n$ divides $b_n$ for each $n\geq2$. Set $r_n=b_n/N_n$ for $n\ge2$.
So
$$\D_{n}\oplus N_n\{0,1,\ldots,r_n\}=\{0,1,\ldots,b_n-1\}$$
for $n\ge2$.
We denote $C_1=\{0\}$ and $C_n= N_n\{0,1,\ldots,r_n\}$ for $n\geq2$ and let $\mu_{{\bf b}, {\bf C}}$ be the Cantor-Moran measure generated by
${\bf b}=\{b_n\}_{n=1}^{\infty}$ and ${\bf C}=\{C_n\}_{n=1}^{\infty}$. Then
 \begin{eqnarray*}
 \mu_{{\bf b},{\bf D}}\ast\mu_{{\bf b}, {\bf C}}&=&\delta_{\frac{1}{b_1}\D_{1}}\ast{\mathcal L}_{[0, \frac1{b_1}]}\\
&=&{\mathcal L}_{[0, \frac{N_1}{b_1}]}.
  \end{eqnarray*}
It implies that the generalized Fuglede's Conjecture holds for $\mu_{{\bf b},{\bf D}}$.

$``{\rm(ii)}\Rightarrow {\rm(i)}"$ Suppose there exists a Borel probability measure $\nu$ such that  $\mu_{{\bf b},{\bf D}}\ast\nu={\mathcal L}_{[0, \frac{N_1}{b_1}]}$.
Let $f(x)=\frac{N_1}{b_1}x$ for $x\in\R$.
Then $\mu_{{\bf b},\frac{b_1}{N_1}{\bf D}}\ast(\nu\circ f)={\mathcal L}_{[0, 1]}$. From Theorem \ref{th-Add1}, $\mu_{{\bf b},\frac{b_1}{N_1}{\bf D}}$ is a spectral measure and so is $\mu_{{\bf b},{\bf D}}$ according to Lemma \ref{lem4}.
\qed

\section{Proof of Theorem \ref{th1.5} and Open Question}
In this section, firstly, we give the proof of Theorem \ref{th1.5}. And then, we conclude this paper with a conjecture on the relationship between the spectral Cantor-Moran measure and the tiling of integers.

We first introduce a theorem from Tijideman \cite{T1995}.

\begin{theo}[\cite{T1995}]\label{th-Add2}
Suppose that $A$ is finite, $0\in A\cap B$ and $A\oplus B=\Z$. If $r$ and $\#A$ are relatively prime, then $rA\oplus B=\Z$.
\end{theo}

\noindent{\bf{ Proof of Theorem \ref{th1.5}}.} It would be easier to prove that ${\bf D}_n=\D_{n}\oplus b_n\D_{n-1}\oplus\cdots\oplus b_2\cdots b_n\D_{1}$ is an integer tile if and only if $N_n$ divides $b_n$ for each $n\geq2$.
Suppose for any $n\geq1$, $b_n=r_nN_n$ for some integer $r_n$, then it is clear that ${\bf D}_n=\D_{n}+b_n\D_{n-1}+\cdots+b_2\cdots b_n\D_{1}$ is a direct sum. Denote $C_1=\{0\}$, $C_n= N_n\{0, 1, \ldots, r_n\}$ and
$${\bf C}_n=C_{n}\oplus b_nC_{{n-1}}\oplus\cdots\oplus b_2\cdots b_nC_{1}$$
for $n\geq2$.
Then $\D_{n}\oplus C_n=\{0, 1, \ldots, b_n-1\}$ and so
$${\bf D}_n\oplus{\bf C}_n=\{0, 1, \ldots, N_1b_2\cdots b_n-1\}$$
for each $n\ge1$. It implies that ${\bf D}_n$ is an integer tile.

Suppose that the converse conclusion is false, ${\bf D}_n=\D_{n} + b_n\D_{n-1}+\cdots+ b_2\cdots b_n\D_{1}$ is a direct sum and an integer tile, but there is $N_n$ does not divide $b_n$. Then there must be a prime factor $p_n$ of $N_n$ which is co-prime with $b_n$. Note that
\begin{equation}\label{eq-L3}
  \D_{n}=\{0, 1, \ldots, p_n-1\}\oplus p_n\{0, 1, \ldots, N_n/p_n-1\}.
\end{equation}
The integer tiling property of ${\bf D}_n$ implies that we can find $0\in C\subset\Z$ such that
$${\bf D}_n\oplus C=\{0, 1, \ldots, p_n-1\}\oplus b_n{\bf D}_{n-1}\oplus p_n\{0, 1, \ldots, N_n/p_n-1\}\oplus C=\Z.$$
Take $A=\{0, 1, \cdots, p_n-1\}$ and $B=b_n{\bf D}_{n-1}\oplus p_n\{0, 1, \cdots, N_n/p_n-1\}\oplus C$. As $\gcd(p_n,b_n)=1$ and $\#A=p_n$, from Theorem \ref{th-Add2} we have
$$b_nA\oplus b_n{\bf D}_{n-1}\oplus p_n\{0, 1, \cdots, N_n/p_n-1\}\oplus C =\Z.$$
The above direct sum cannot happen as $\{0, 1\}\subset A\cap {\bf D}_{n-1}$. We have a contradiction and this shows our desired statement holds.
\qed

%
\medskip

Theorem \ref{th1.5} shows that there is a strong link between the spectral Cantor-Moran measure and the integer tile.
For self-similar case ($(b_n, \D_n)\equiv (b, \D)$), ${\L}$aba and Wang \cite{LW2002} conjectured that if $\mu_{b, \D}$ is a spectral measure, then $\D$ is an integer tile. However, integer tile is not a sufficient condition. A simple example observed by  the first-named author, He and Lai\cite{AHL2015}: let $b=8$ and $\D=\{0, 1, 8, 9\}$. Then $\D$ tiles $\Z_{16}$, but the self-similar measure $\mu_{8, \D}$ is not spectral. It mainly because of $\D+8\D$ is not a direct sum.
Recall that $\mu_{\textbf{B},\textbf{D}}=\mu_{n}\ast\mu_{>n}$ and
$$
\mu_n=\delta_{\frac{1}{b_1}\D_{1}}\ast\delta_{\frac{1}{b_1b_2}\D_{2}}\ast\cdots\ast\delta_{\frac{1}{b_1b_2\cdots b_n}\D_{n}}=\delta_{\frac{1}{b_1b_2\cdots b_n}{\bf D}_n},
$$
where ${\bf D}_n=\D_{n}+b_n\D_{{n-1}}+b_2\cdots b_n\D_{1}$.
According to Fuglede's conjecture on cyclic groups, one may ask:

  {\bf Qu:} If $\mu_{\textbf{B},\textbf{D}}$ is a spectral measure, is each iterated digit set ${\bf D}_n=\D_{n}+b_n\D_{{n-1}}+b_2\cdots b_n\D_{1}$ $(n\geq1)$ an integer tile?

 Our Theorem \ref{th1.5} settles the case where $\D_n=\{0, 1, \ldots, N_n-1\}$. And we look forward to more positive results regarding this question.

\end{document}